\documentclass{cimart}

\newtheorem*{questionB}{Problem B}
\newtheorem*{questionD}{Theorem of Stäckel}
\newtheorem*{questionBD}{Problem B involving derivatives}
\newtheorem*{TR}{Rouché's Theorem}
\newtheorem*{Teorema de H}{Hurwitz's Theorem}

\newcommand{\A}{\mathbb{A}}
\newcommand{\C}{\mathbb{C}}
\newcommand{\Q}{\mathbb{Q}}
\newcommand{\QQ}{\overline{\Q}}

\title{On a version of the problem B of Mahler involving derivatives}

\author{Diego Alves}

\authorinfo{Departamento de Ensino, Instituto Federal do Ceará, Brazil}{diego.costa@ifce.edu.br}

\abstract{In 1902, Paul Stäckel constructed an analytic function $f(z)$ in a neighborhood of the origin, which was transcendental, and with the property that both $f(z)$ and its inverse, as well as its derivatives, assumed algebraic values at all algebraic points in this neighborhood. Inspired by this result, Mahler in 1976 questioned the existence of a transcendental entire function $f(z)$ such that $f(\QQ)$ and $f^{-1}(\QQ)$ are subsets of $\QQ.$ This problem was solved by Marques and Moreira in 2017. As St\"acklel's result involved derivatives, it is natural to question whether we have an analogous result for transcendental entire functions involving derivatives. In this article, we show that there are an uncountable amount of such functions.}



\VOLUME{32}
\YEAR{2024}
\ISSUE{1}
\firstpage{305}
\DOI{https://doi.org/10.46298/cm.12784}

\begin{document}

\section{Introduction}
       
An analytic function $f$ over a domain $\Omega\subseteq \C$ is said to be a \textit{transcendental function} over $\C(z)$, if the only polynomial complex $P$ in two variables satisfying $P(z,f(z))=0,$ for all  $z \in \Omega,$ is the null polynomial. From the well-known characterization that an entire function is transcendental if and only if it is not a polynomial function, we obtain several examples of transcendental functions, such as the sine and cosine trigonometric functions and the exponential function. 

At the end of the 19th century, after the proof of the Hermite-Lindemann Theorem, which demonstrates that the exponential function assumes transcendental values in every non-zero algebraic number, the question arose: does a transcendental entire function almost always (e.g. with the exception of a finite subset of $\QQ$) assume transcendental values in algebraic points?

In this sense, Straüs in 1886, tried to prove that an analytic and transcendental function in a neighborhood of the origin cannot assume rational values in all rational points of its domain. However, he was surprised when Weierstrass presented him with a counterexample in the same year. Weierstrass also conjectured that there were transcendental entire functions that assumed algebraic values at all algebraic points. This conjecture was confirmed by Stäckel \cite{stackel1895} who proved the following more general result: if $\Sigma$ is a countable subset of $\C$ and $T$ is a dense subset of $\C,$ then there exists an entire and transcendental function $f$ such that $f(\Sigma) \subseteq T.$ The Weierstrass conjecture is obtained by taking $\Sigma=T= \QQ.$ In this context of arithmetic behavior of analytic functions, Stäckel has other results, of which we highlight the following, proved in \cite{stackel1902}
\begin{questionD} There exists a transcendental function $$f(z)=-z+\displaystyle\sum_{h=2}^{\infty}f_h z^h,$$ with rational coefficients, which converges in a neighbourhood of the origin and has the property that both $f(z)$ and its inverse function, as well as all their derivatives, are algebraic at all algebraic points in this neighbourhood.
\end{questionD}

Based on this theorem, Mahler \cite{mahler1976} questioned whether this result (without involving derivatives) could be extended to transcendental entire functions, and named this question Problem B. Precisely, we transcribe below the statement of the question proposed by Mahler:

\begin{questionB}
          Does there exist an entire transcendental function $$f(z)= \sum_{h=0}^{\infty} f_h z^h$$ with rational  coefficients $f_h$ such that both $f(z)$ and its inverse function are algebraic in all algebraic points? 
         \end{questionB}
         Since, by Picard's Great Theorem, an entire and transcendental function cannot even be one-to-one (let alone admit an inverse), it is believed that Mahler when enunciating the Problem $B$, was actually referring to the inverse image instead of the function inverse, as it might seem at first reading. This understanding was adopted by Marques and Moreira \cite{Marques2017}, who in 2017 solved the Problem B of Mahler proving the following result:
          \begin{theorem}[Cf. Theorem 1 of \cite{Marques2017}]
         There are uncountably many transcendental entire functions $f(z)$ with rational coefficients such that $$f(\overline{\mathbb{Q}}) = f^{-1}(\overline{\mathbb{Q}}) = \overline{\mathbb{Q}},$$ where $f^{-1}(\QQ)$ is the inverse image of $\QQ$ by $f.$
          \end{theorem}
          Since Stäckel's result involved derivatives, it is natural to ask whether we could extend the previous theorem to a version involving derivatives, as follows:
        \begin{questionBD}
          Does there exist an entire transcendental function $f(z)$ with rational  coefficients such that
         $$\left[f^{(j)}(\overline{\mathbb{Q}}) \cup \left(f^{(j)}\right)^{-1}(\overline{\mathbb{Q}})\right] \subseteq \overline{\mathbb{Q}}, \, \, \text{for all} \, \,  j \geq 0?$$  
\end{questionBD}

Here $f^{(j)}$ denotes the derivative of order $j$ of $f$ with $f^{(0)}=f$ and $\left( f^{(j)}\right)^ {-1}(\QQ)$ denotes the inverse image of $\QQ$ by the function $f^{(j)}.$

We refer the reader to \cite{mahler1976,waldschmidt2003} (and references therein) for more about this subject.

In this paper, inspired by the ideas of Marques and Moreira \cite{Marques2017}, and using the theorems of Rouché and Hurwitz, we provide an affirmative answer to the question above, proving the following theorem:
 \begin{theorem}\label{aqf}
          There are uncountably many transcendental entire functions $f(z)$ with rational coefficients such that $$f^{(j)}(\overline{\mathbb{Q}}) \subseteq \overline{\mathbb{Q}} \, \,  \text{and} \, \, \left(f^{(j)}\right)^{-1}(\overline{\mathbb{Q}}) \subseteq \overline{\mathbb{Q}}, \, \, \text{for all} \, \,  j \geq 0.$$
          \end{theorem}
        For a more didactic exposition and for the sake of completeness, we will review some basic facts about zeros of analytic functions in a complex variable. Under certain assumptions, it is possible to compare the number of zeros of two (or more) analytic functions in a domain in which both are defined. It is on this topic that the next two theorems address.
        \begin{TR}
        Let $U \subseteq \mathbb{C}$ be a domain and suppose that $f,g:U \rightarrow \mathbb{C}$ are analytic functions. Let $V \subseteq U$ be a bounded and closed domain whose boundary $\partial V$ is a piece-wise smooth Jordan curve such that $V - \partial V$ is a domain. If \begin{equation}\label{kk}
         |f(z)-g(z)| < |f(z)|, \, \, \text{for all} \, \, z \in \partial V, 
         \end{equation}
         then $f$ and $g$ have the same number of zeros inside $V,$ each of them counted as many times as their multiplicity.
         \end{TR}
         \proof See \cite[Chapter 6]{Soares}.
      
         \begin{Teorema de H}
          Let $G$ be a domain and suppose that $\{f_n:G \rightarrow \C\}_n$ is a sequence of analytic functions on $G$ that uniformly converge on compact subsets of $G$ to $f.$ If $f \not\equiv 0$ and there is  $R>0$ such that $\overline{B}(a,R) \subseteq G$ and furthermore $f$ does not vanish at any point of $\partial B(a,r),$ then there is an integer $N$ such that $f$ and $f_n$ have the same number of zeros in $B(a;R).$ 
           \end{Teorema de H}
           \proof See  \cite[Chapter VII]{Conway}.

           \medskip
           
         In the applications made in this work, we will work with integer functions and the set V will always be a closed ball. In Rouché's Theorem notations, we will have $U = \C$ and $V=\overline{B}(a,R),$ therefore, to apply Rouché's Theorem, our only concern will be to check whether the inequality (\ref{kk}) is satisfied.
         
        Now, we are ready to demonstrate Theorem \ref{aqf}.

         \section{The proof}
        To simplify our exposition, we will denote the set of real algebraic numbers by $\A,$ and given two integers $a$ and $b$ with $a<b,$ [a,b] to denote the set $\{a,a+1, \ldots, b\}.$ Let us consider $\{ \alpha_1, \alpha_2, \alpha_3, \ldots \}$ an enumeration of $\QQ$ such that $\alpha_1=0$ and for every $n \geq 1,$
         \begin{enumerate}  
        \item[$\bullet$] $\alpha_{3n-1}$ and $\alpha_{3n}$ are non-real complex numbers with $\alpha_{3n-1}= \overline{\alpha}_{3n};$ 
         \item[$\bullet$] $\alpha_{3n+1}$ is a real number.
         \end{enumerate}
         We will construct the desired function inductively. We start with $f_1(z) = z^2 $ and $P_1(z)=1.$ Note that $f_1(\{\alpha_1\})=f_1^{-1}(\{\alpha_1\})=\{0\} \subseteq \overline{\mathbb{Q}}$. We will construct a sequence of analytic functions $f_2(z), f_3(z), \ldots$ of the form $$ f_m(z)= f_{m-1}(z) + \epsilon_m z^{m+1} P_m(z) = \displaystyle\sum_{i=2}^{t_m} a_i z^i,$$ satisfying
         \begin{enumerate}
         \item[(i)] $f_{m}(z), P_m (z) \in \A[z]-\{0\}$ and  $a_{t_m} \neq 0,$ where $t_m \geq m + 1;$
         \item[(ii)] $P_{m-1}(z) | P_m(z)$ and $P_m(0) \neq 0;$
        \item[(iii)] $\epsilon_m \in \A;$
        \item[(iv)] $0 < |\epsilon_m| < \dfrac{1}{L(P_m)m^{m+1+\deg P_m}}:= \Gamma_m;$
        \item[(v)] $a_{2}, \ldots, a_{m+1} \in \mathbb{Q} - \{0\},$
        \end{enumerate}
        where $L(P_m)$ denotes the length of the polynomial $P_m,$ given by the sum of absolute values of its coefficients. The desired function will have the form $$f(z) = z^2 + \sum_{n \geq 2} \epsilon_n z^{n+1} P_n(z).$$
        Now using that $|P(z)| \leq L(P) \max \{1, |z|\}^{\deg P},$ for all $P(z) \in \mathbb{C}[z],$ we get that for all $z \in B(0,R)$ (where R is a fixed arbitrary positive real number) $$|\epsilon_n z^{n+1} P_n(z)| <  \left( \dfrac{\max \{1,R\}}{n} \right)^{n+1+ \deg P_n} :=u_n.$$
        Since $u_n < (1/2)^n,$ for all sufficiently large $n,$ it follows by the comparison test and by the Weierstrass M-test that the series $$f(z) = z^{2} + \displaystyle\sum_{n \geq 2} \epsilon_n z^{n+1} P_n(z)$$ converges uniformly on each ball $B(0,R) $ (which implies uniform convergence on compact subsets of $\mathbb{C}).$ In particular, $f$ is an entire function.
        
        Suppose we have constructed $f_n$ satisfying conditions (i) to (v). Now, we will construct $f_{n+1}$ with the desired properties. Since $f_n$ is a polynomial of degree greater than $n,$ the set $A_n$ defined below is finite
        $$A_n = \displaystyle\bigcup_{j=0}^{n} (f_n^{(j)})^{-1}(\{ \alpha_1, \ldots, \alpha_{3n+1}\}) = \{0, y_1, \ldots, y_s\}.$$
        So, we define 
         $$f_{n+1}(z) = f_n(z) + \epsilon_{n+1}z^{n+2}P_{n+1}(z),$$
         where $$P_{n+1}(z) = P_n(z) (z - \alpha_2) \cdots (z-\alpha_{3n+1}) \displaystyle\prod_{i=1}^{s} (z-y_i)^{\deg f_n+n+1}$$ and $\epsilon_{n+1}$ will be chosen later. Note that $P_{n+1}(z) \in \A[z].$ Indeed, if $y_k$ is not real, as $y_k$ is a zero of $f_n^{(j)}(z) - \alpha_i,$ for some $j \in [0, \ldots,n]$ and $i \in [1, 3n+1]$ and $f_n(z) \in \A[z],$ we have that $\overline{y}_k$ is a zero of $f_n ^{(j)}(z) - \overline{\alpha}_i.$ Since $\{\alpha_1, \ldots, \alpha_{3n+1}\}$ is closed for complex conjugation, it follows that $\overline{y}_k \in A_n.$ Thus, $P_{n+1}(z) \in \A[z]$ since $P_n(z) \in \A[z]$ and the set of algebraic numbers $\{\alpha_2, \ldots, \alpha_{3n+1}, y_1, \ldots, y_s\}$ is composed of real numbers and pairs of conjugate of complex (non-real) numbers. Also $P_{n+1}(0) \neq 0,$ since $P_n(0) \neq 0$ and $\alpha_2,\ldots, \alpha_{3n+1},y_1,\ldots, y_s$ are not null. As $A_n$ is a finite set, we can choose a positive real number $r_{n+1}$ such that $$n+1 < r_{n+1}< n+2 \, \, \text{and} \, \, A_n \cap \partial B(0,r_{n+1}) = \emptyset.$$ 
         So, for all $j \in [0,n]$ and $i \in [1,3n+1],$ we have $\displaystyle\min_{|z|= r_{n+1}} |f^{ (j)}_n(z)-\alpha_i|>0.$ In this way, we can choose $\epsilon_{n+1}$ satisfying
         \begin{equation}\label{pcon}
        |\epsilon_{n+1}| < \dfrac{\displaystyle\min_{|z|= r_{n+1}} |f^{(j)}_n(z)-\alpha_i|}{\displaystyle\max_{0 \leq k \leq n} \displaystyle\max_{|z|=r_{n+1}} |\left( z^{n+2}P_{n+1}(z) \right)^{(k)}|} = \Lambda_{i,j},
        \end{equation}
         for all $i \in [1,3n+1]$ and for all $j \in [0,n].$ Therefore, for any $j \in [0,n]$ and $i \in [1,3n+1],$ we obtain that for all $w \in \partial B(0,r_{n+1})$ 
         \begin{eqnarray}
        \left| f_n^{(j)}(w) - \alpha_i \right|& \geq &  \displaystyle\min_{|z|= r_{n+1}} |f^{(j)}_n(z)-\alpha_i|  \nonumber\\
        & > & |\epsilon_{n+1}| |\displaystyle\max_{0 \leq k \leq r} \displaystyle\max_{|z|=r_{n+1}} |\left( z^{n+2}P_{n+1} (z)\right)^{(k)}| \nonumber\\
        &\geq& |\epsilon_{n+1}|  \left|\left[\left( z^{n+2}P_{n+1}(z) \right)^{(j)} \right]_{|z=w} \right|\nonumber\\
        &=& |f_{n+1}^{(j)}(w) - \alpha_i - ( f_n^{(j)}(w) - \alpha_i) |. \nonumber
        \end{eqnarray}
        Then, by Rouché's theorem, $f^{(j)}_n(z) - \alpha_i$ and $f^{(j)}_{n+1}(z) - \alpha_i$ have the same number of zeros (counting multiplicities) in $B(0,r_{n+1}),$ for any $j \in [0,n]$ and $i \in [1,3n+1].$  Note that, if $0$ is a root of multiplicity $m \geq 1$ of $f^{(j)}_n(z) - \alpha_i,$ then $0$ is a root of multiplicity $m$ of $f^ {(j)}_{n+1}(z) - \alpha_i.$ In fact, if $j=0,$ since $f_n^{(0)}(z)=f_n(z)$ and $f_{n+1}^{(0)}(z)=f_{n+ 1}(z)$ has the same coefficients up to order $2,$ with the coefficient of order $2$ not zero, it follows that $m=2$ and $0$ is a multiplicity root $m=2$ of $f_{ n+1}^{(j)}(z) - \alpha_i.$ If $j \geq 1,$ then the multiplicity $m$ is equal to $1,$ since $$ \left[\left(f_n^{(j)}(z) - \alpha_i \right)' \right]_{|z=0} = f_n^{(j+1)}(0) =(j+1)! a_{j+1},$$ which is not zero because $2 \leq j +1 \leq n+1$ and $a_2, \ldots, a_{n+1}$ are non-zero. Since the coefficients of $f_n(z)$ and $f_{n+1}(z)$ coincide up to the order $n+1,$ we have $f_n^{(j)}(0) = f_{n +1}^{(j)}(0),$ for all $j \in [0,n+1],$ so that $f_{n+1}^{(j)}(0) - \alpha_i = f_{n}^{(j)}(0)- \alpha_i =0$ and $$\left[\left(f_{n+1}^{(j)}(z) - \alpha_i \right)' \right]_{|z=0} = f_{n+1}^{(j+1)} (0) =  f_{n}^{(j+1)}(0) = (j+1)!a_{j+1} \neq 0.$$  Therefore, $0$ is a root of multiplicity $m =1$ of $f^{(j)}_{n+1}(z) - \alpha_i.$ On the other hand, if an element $\lambda \in B(0,r_{n+1})-\{0\}$ is a zero of multiplicity $m \geq 1$ of $f_{n}^{(j) }(z) - \alpha_i,$ then $\lambda = y_l$ for some $l \in \{1, \ldots,s\},$ so that $\lambda$ is a zero of multiplicity greater than $\deg f_n$ of $\left[ z^{n+2} P_{n+1}(z)\right]^{(j)},$ since $j \leq n.$ Therefore, $\lambda$ is also a zero of multiplicity $m$ of $f^{(j)}_{n+1}(z)- \alpha_i.$
        From these facts, we obtain that the polynomials $f^{(j)}_{n+1}(z)- \alpha_i$ and $f^{(j)}_{n}(z)- \alpha_i$ (for any $j \in [0,n]$ and $i \in [1, 3n+1]$) have exactly the same zeros with the respective multiplicities in $B(0,r_{n+1}).$ In particular, for all $i \in [1,3n+1]$ and $j \in [0,n],$ we have 
        \begin{equation}\label{pp}
        (f_n^{(j)})^{-1} ( \alpha_i) \cap B(0,r_{n+1}) = (f_{n+1}^{(j)})^{-1} ( \alpha_i) \cap B(0,r_{n+1}).
        \end{equation}
        This argument ensures that no new pre-image from the set $\{\alpha_1, \ldots, \alpha_{3n+1}\}$ by $f^{(j)}_{n+1}$ in $B(0,r_{n+1}),$ with $j \in [0,n],$ will appear in addition to those already existing by the function $f^{(j)}_n.$ Note also that, since $f_{n+1} \in \A[z]$ has degree greater than $n$ and $\overline{\mathbb{Q}}$ is algebraically closed, $f^{(j) }_{n+1} ( \{\alpha_1, \ldots,\alpha_{3n+1}\})$ and $\left( f^{(j)}_{n+1} \right)^{ -1}( \{\alpha_1, \ldots,\alpha_{3n+1}\}),$ with $j \in [0,n],$ are subsets of $\overline{\mathbb{Q}}.$
        
        Let's write $$ f_{n+1}(z) = \sum_{i=2}^{t_{n+1}} a_i z^i.$$

        The coefficients of $f_{n+1}$ from the order $2$ to the order $n+1$ coincide with the coefficients of $f_n,$ and therefore are non-zero rational numbers. We show that it is possible to choose $\epsilon_{n+1} \in \A$ satisfying (\ref{pcon}), (iv) and such that $a_{n+2}$ is also a non-zero rational number. In fact, let $c_{n+2}$ be the coefficient of $z^{n+2}$ in $f_n(z),$ then $$a_{n+2} = c_{n+2} + \epsilon_{n+1}P_{n+1}(0).$$ Since $P_{n+1}(0) \neq 0,$ we can choose $p/q \in \mathbb{Q} - \{ 0\}$ such that
        $$0< \left|c_{n+2}- \frac{p}{q}\right| < |P_{n+1}(0)| \min \{ \Gamma_{n+1}, \Lambda_{1,0}, \ldots, \Lambda_{3n+1,0}, \ldots, \Lambda_{1,n}\ldots, \Lambda_{3n+1,n} \}.$$
        We define $\epsilon_{n+1} = \left(p/q - c_{n+2}\right)/P_{n+1}(0),$ then $\epsilon_{n+1} \in \A,$ satisfies (iv) and (\ref{pcon}) and furthermore, $a_{n+2} = p/q \in \mathbb{Q} - \{0\}.$ So, by construction, the function $$f(z) = z^{2} + \displaystyle\sum_{n = 2}^{\infty} \epsilon_{n} z^{n + 1} P_n(z) = \lim_{n \to \infty} f_n(z)$$ is an entire function and has rational coefficients. Now, let's check that
        \begin{equation}\label{abc}
        f^{(j)}( \overline{\mathbb{Q}}) \subseteq \overline{\mathbb{Q}} \, \, \text{and} \, \, (f^{(j)})^{-1} (\overline{\mathbb{Q}}) \subseteq \overline{\mathbb{Q}}, \, \, \, \text{for all} \, \, j \geq 0.
         \end{equation}
         To see the first inclusion, notice that by construction
         $$f_n( \alpha_i) = f_i(\alpha_i), \, \, \text{for all} \, \, n \geq i.$$
         In this way, we have $$f(\alpha_i) = \lim_{n \to \infty} f_n(\alpha_i)=f_i(\alpha_i) \in \overline{\mathbb{Q}}.$$
         Also, as the multiplicity of $\alpha_i$ in $P$ tends to infinity when $n$ tends to infinity, we have (by the general Leibniz rule for differentiation) a similar property involving the derivatives of $f$ and of the functions $f_n$, precisely: given a positive integer $j$, there exists $m= m(j,i) \in \mathbb{N}$ such that $$f_n^{(j)}(\alpha_i) = f_m^{(j)}(\alpha_i) \, \, \text{for all} \, \, n \geq m.$$
         As $f_n$ converges uniformly on compact subsets of $\C$ for f the same kind of convergence occurs from $f_n^{(j)}$ to $f^{(j)},$ for all $j \geq 1,$ so that $$f^{(j)}(\alpha_i) = \lim f^{(j)}_n (\alpha_i)= f^{(j)}_{m(j,i)}(\alpha_i) \in \overline{\mathbb{Q}},$$ showing the first inclusion in (\ref{abc}).

        To prove the second inclusion in $(\ref{abc}),$ we start by noting that if $i \leq 3n+1$ and $t,j \leq n,$ then
         \begin{equation}\label{xx}
        (f_n^{(j)})^{-1} ( \alpha_i) \cap B(0,r_t) = (f_{n+1}^{(j)})^{-1} ( \alpha_i) \cap B(0,r_t).
        \end{equation}
        This follows directly from (\ref{pp}) and from $r_t \leq r_{n+1}.$ Therefore, if $i,t$ are positive integers and $j$ is a non-negative integer and we set $l= \max\{t,i,j\},$ then we have
         \begin{equation}\label{sx}
             (f^{(j)}_n)^{-1}(\alpha_i) \cap B(0, r_t) = (f^{(j)}_l)^{-1}(\alpha_i) \cap B(0,r_t) \, \, \text{for all} \, \, n \geq l.
         \end{equation} 
        Indeed, if $n=l$, then the result is obvious.  Hence, we can suppose $n \geq l+1.$ In this case, $\max\{t,i,j\} = l \leq n-1,$ consequently $i \leq 3l+1 \leq 3(n-1)+1,$ so by (\ref {xx}),
         \begin{eqnarray}
        \left(f_l^{(j)}\right)^{-1}(\alpha_i) \cap B(0,r_t) & = &   \left(f_{l+1}^{(j)}\right)^{-1}(\alpha_i) \cap B(0 , r_t) \nonumber\\
         &  &  \, \, \, \, \, \, \, \, \, \, \, \, \, \, \, \, \, \, \, \, \, \, \, \, \, \, \, \, \, \, \, \, \vdots \nonumber\\
        & = & \left(f_{n-1}^{(j)}\right)^{-1}(\alpha_i) \cap B(0 , r_t), \nonumber \\
        & = & \left(f_{n}^{(j)}\right)^{-1}(\alpha_i) \cap B(0, r_t), \nonumber
        \end{eqnarray}
        proving what we had asserted. From this it follows that
        $$(f_l^{(j)})^{-1}(\alpha_i) \cap B(0,r_t) \subseteq (f^{(j)})^{-1}(\alpha_i) \cap B(0, r_t), \, \, \text{where} \, \, l = \max\{t,i,j\}.$$
        In fact, let $y \in (f_l^{(j)})^{-1}(\alpha_i) \cap B(0, r_t),$ then we have $f_n^{(j )}(y)=\alpha_i$ for all $n \geq l,$ and consequently, $f^{(j)}(y) = \displaystyle\lim_{n \to \infty} f^{(j) }_n(y) = \alpha_i.$ Therefore, $y$ belongs to $(f^{(j)})^{-1}(\alpha_i) \cap B(0,r_t).$ We claim that 
        \begin{equation}\label{zz}
         (f_l^{(j)})^{-1}(\alpha_i) \cap B(0,r_t) = (f^{(j)})^{-1}(\alpha_i) \cap B(0,r_t), \, \, \text{where} \, \, l = \max\{t,i,j\}.
        \end{equation}

       Let us suppose by absurdity that there is 
       \[
            w \in (f^{(j)})^{-1}(\alpha_i) \cap B(0 , r_t)\setminus (f_l^{(j)})^{-1 }(\alpha_i) \cap B(0 , r_t),
        \]
        then $w$ would be at a positive distance $\delta$ from the finite set $(f_l^{(j)})^{-1}(\alpha_i) \cap B(0,r_t).$ We set $$  S = (f^{(j)})^{-1}(\alpha_i) \bigcup \left[ \displaystyle\bigcup_{n \geq l} (f_n^{(j)})^{-1}(\alpha_i) \right].$$ Since $f^{(j)}$ is a non-constant integer and for each $n \geq l,$ $f_n$ is a polynomial of degree greater than $j,$ we have that $S$ is countable. Therefore, it is possible to choose $\beta < \delta$ such that $$B= B(w, \beta) \subseteq B(0,r_t) \, \, \text{and} \, \, S \cap \partial B = \emptyset.$$ Since $f_n^{(j)}(z) - \alpha_i$ converges uniformly on compact subsets of $\mathbb{C}$ for the entire function $f^{(j)}(z) - \alpha_i,$ which does not vanish at any point of $\partial B,$ follows by Hurwitz's Theorem that $f^{(j)}(z)- \alpha_i$ and $f_n^{(j)}(z) - \alpha_i$ has the same number of zeros in $B(w, \beta)$ for all sufficiently large $n.$ On the other hand, by (\ref{sx}), the zeros of $f^{(j)}_n(z) - \alpha_i$ and $f_l^{(j)}(z) - \alpha_i$ in $B(0,r_t )$ coincide for $n \geq l.$ So, since $f_l^{(j)}(z) - \alpha_i$ has no zeros $B(w,\beta),$ it follows that $f^{(j)}(z)-\alpha_i$ does not have zeros in $B(w,\beta).$ However, $w \in B(w, \beta)$ and $f^{(j)}(w)-\alpha_i = 0,$ a contradiction. Therefore the equality (\ref{zz}) holds. From this it follows that $(f^{(j)})^{-1} (\overline{\mathbb{Q}}) \subseteq \overline{\mathbb{Q}},$ for all $j \geq 0.$ In fact, using (\ref{zz}), we obtain
        \begin{eqnarray}
        (f^{(j)})^{-1}(\alpha_i) & = &   \bigcup_{t=\max\{i,j\}}^{\infty} (f^{(j)})^{-1}(\alpha_i) \cap B(0,r_t)  \nonumber\\
         & = & \bigcup_{t=\max \{i,j\}}^{\infty} (f_t^{(j)})^{-1}(\alpha_i) \cap B(0, r_t) \nonumber
        \end{eqnarray}
        which is a subset of $\QQ.$

       \medskip
       
        By construction, $f$ is an entire function and is not a polynomial (because its coefficients of order greater than to $1$ are non-zero), so $f$ is a transcendental function. Finally, note that there is an uncountable amount of possible choices for $f,$ because at each step, we have infinite possibilities of choices for $\epsilon_{n+1}$ and thus also for $a_{n+2}.$ In this way, we construct an uncountable amount of transcendental entire functions $f$, with rational coefficients and such that $$\left[ f^{(j)}(\overline{\mathbb{Q}}) \cup  (f^{(j)})^{-1}(\overline{\mathbb{Q}}) \right] \subseteq \overline{\mathbb{Q}}, \, \, \text{for all} \, \, j \geq 0.$$     
         \qed

%%% REFERENCES %%%

\EditInfo{January 5, 2024}{February 10, 2024}{Lenny Fukshansky}

\end{document}